# FFT-based Computation of Polynomial Coefficients and Related Tasks

Author address    Rollenbergstrasse 2, CH-8463 Benken ZH, hr@thomannconsulting.ch
Version          12.07.2016 10:22

**CONTENTS**

**Abstract**

We present a FFT-based algorithm for the computation of a polynomial's coefficients from its roots, and apply it to obtain the coefficients of interpolation polynomials, to invert Vandermondians and to evaluate the symmetric functions of a set of parameters. Analytic and numerical evidence with problem sizes up to and beyond $n = 2000$ confirms that it is superior over previous algorithms for these problems in case of parameters taken from uniform or almost uniform distributions on or near circles in the complex plane, and of comparable performance in other cases. Its time complexity is $O(n^2)$ and its storage complexity for tasks other than Vandermondian inversion $O(n)$.

# 1 Introduction

The topic of this paper is a new, FFT-based algorithm for the computation of the coefficients of a polynomial from its roots, and its application to the computation of the coefficients of interpolation polynomials, the inversion of Vandermondians and the symmetric functions of a set of parameters.

Computation of the polynomial coefficients from the roots is indeed the critical part of common algorithms for the determination of interpolation polynomials as well as the inversion of Vandermondians. The polynomial coefficients are (up to sign) the symmetric functions of the roots. Therefore accurate computation of the polynomial coefficients enhances the numerical performance of these tasks.

We start by introducing some problems and related algorithms. Clause 2 defines the evaluation methodology and discusses the conditions of these problems. Clause 3 contains the evaluation and comparison of the numerical performance of our new versus the conventional solution of the main problem, using polynomials whose coefficients are given by closed formulae. In clause 4 polynomials with more general roots are treated, and the performance is measured against interpolation. In clause 5 interpolation at the roots is evaluated, which is directly related to the performance for the inversion of Vandermondians. Clause 6 evaluates the application of our method to the calculation of the coefficients of interpolation polynomials. The findings are briefly summarized in the last clause.

## 1.1 Main problem and algorithm: Polynomial coefficients

We denote by problem A the computation of the coefficients of the $n^{th}$-order polynomial

(1) $P(\omega) = \prod \omega - \omega_k = \sum a_m \omega^m$

from its roots $\omega_1, \cdots, \omega_n$. Throughout this paper we require all roots to be non-zero. Thus $a_0$ and $a_n$ are non-zero.

If $\omega^N = 1$ $(N > n)$, then the Fourier Transform yields

(2) $a_m = \frac{1}{N} \sum_{0 \leq j < N} P(\omega^{-j}) \omega^{jm}$ $(0 \leq m \leq n)$.

Indeed, the r.h.s. equals $\frac{1}{N} \sum \sum a_k \omega^{-jk} \omega^{jm} = a_k \delta_{m-k}$ $(0 \leq m \leq n)$. It vanishes for $n < m < N$.

This gives rise to

    <u>Algorithm P:</u>    FFT-based solution of problem A

    Input:    Complex roots $\omega_1, \cdots, \omega_n$

    Output: $a_m$ $(0 \leq j \leq n)$

    Computation:

        a) Choose $N: n < N \leq 2n$
        b) For $j = 0, 1, \cdots, N-1$ do $p_j = \prod (\omega^{-j} - \omega_k)$
        c) $a_0, \cdots, a_n 0_{n+1}, \cdots, 0_{N-1} = FFT(p_0, \cdots, p_{N-1})$
        d) Output $a_0 \cdots a_n$

The decimation to be used in the FFT determines the choice of N and in turn the time and storage complexity. In our implementation we stick to radix-2 decimation and use the simple Cooley-Tukey procedure, while the sophisticated FFT implementation of Mathematica allows $N = n + 1$. Algorithm P requires $nN + O(N \log N)$ complex arithmetic additions and multiplications and $N + O(1)$ complex storage.

The leading terms amount to $8n^2 < 8nN \leq 16n^2$ of floating point operations and $n < N \leq 2n$ floating point storage. In case of real roots the symmetry $p_{N-j} = \bar{p}_j$ allows to half the complexity to $4n^2 < 4nN \leq 8n^2$ floating point operations and $n/2 < N/2 \leq n$ floating point storage.

Unlike the next one, our algorithm is insensitive to the ordering of the roots. Indeed, unless some partial products are lying outside the domain of the floating point arithmetic, the values and ordering of the products in step b) are insensitive to and independent from the ordering of the roots.

## 1.2 Alternative algorithms for the main problem

Conventional methods for the computation of polynomial coefficients fall into two categories. The first one is provided by algorithm 1.1 of （Calvetti, 2003）, which is based on the recursion

(3) $a_{j,j} = 1, a_{i,j} = 0\ (0 \leq j \leq n, i > j), a_{m,k} = a_{m,k-1} - a_{m+1,k-1}\omega_k\ (0 \leq m < k \leq n),$

where $\prod_{1 \leq j \leq k} \omega - \omega_k = \sum a_{m,k}\omega^m$. This gives rise to

<u>Algorithm R: Solution of problem A based on recursion</u>

Input:   Complex roots $\omega_1, \cdots, \omega_n$

Output: $a_m\ (0 \leq j \leq n)$

Computation:

 a) $a_n = 1, a_j = 0\ (0 \leq j < n)$
 b) for $k = n-1, n-2, \cdots, 0$ do
   $x = 1$
   for $m = n-1, n-2, \cdots, 0$ do
    $y = a_m$
    $a_m = a_m - \omega_k x$
    $x = y$
   end
  end

This algorithm takes $\binom{n}{2} + O(n)$ complex additions and multiplications and $O(n)$ complex storage. The leading term amounts to $4n^2$ floating point operations, reducing to $n^2$ floating point operations in case of real roots.

While more efficient than algorithm P, it is very sensitive to the ordering of the roots. It's accuracy is much improved when combined with Leja ordering （BaglamaJ, 1998）（Calvetti, 2003）, taking another $O(n^2)$ arithmetic operations.

Leja ordering of the roots $\omega_1, \cdots, \omega_n$ is defined by

(4) $|\omega_1| = \max|\omega_j|, \prod_{i<k}|\omega_i - \omega_k| = \max_{j \geq k} \prod_{i<k}|\omega_i - \omega_j|.$

The following algorithm establishes the Leja order:

<u>Algorithm L: Leja ordering of a set of values</u>

Input:   $\omega_1, \cdots, \omega_n$

Ouput:   Parameter set re-arranged in Leja order.

Computation:

 a) Set $p_i = 1\ (i = 1, \cdots, n)$
 b) Find $k$: $|\omega_k| = \max|\omega_j|$
 c) Exchange $\omega_k$ and $\omega_1$
 d) for $k = 2, \cdots, n-1$ do
   for $j = k, \cdots, n$ do
    $p_j = p_j * |\omega_{k-1} - \omega_j|^2$
   end
   Find $i$: $|p_i| = \max_{j \geq k}|p_j|$
   Exchange $\omega_k$ and $\omega_i$
   Exchange $p_k$ and $p_i$
  end

The square absolute value is used for complex roots because it avoids taking the square root. For real roots the absolute value is the suitable choice.

Algorithm L takes $\binom{n}{2}$ complex additions, square absolute value evaluations, real multiplications and maximum function evaluations, and $O(n)$ real storage. The leading term amounts to $\frac{7}{2}n^2$ floating point operations, reducing to $2n^2$ floating point operations in case of real roots.

We denote by algorithm R+ the combination of algorithm L followed by algorithm R. The leading term of the complexity amounts to $\frac{15}{2}n^2$ floating point operations for complex roots, $3n^2$ floating point operations for real roots.

When algorithm P is optimally implemented, then on complex roots algorithm R+ takes $\frac{1}{16}$ less time, on real roots $\frac{1}{4}$ less time. The numerical performance of the two algorithms is subject of the subsequent clauses.

Algorithms of the second category solve a linear system with Vandermondian coefficient matrix. When using a standard solution algorithm they have $O(n^3)$ time complexity. The Björck-Pereyra algorithm （A. Björck, 1970）takes only $O(n^2)$ time. However, in the special case of data at the roots plus another point, such that the constant vector is a multiple of the vector $(1,0,\cdots,0)$, this algorithm reduces to algorithm R. The reader will easily verify this studying the pseudo-code in clause 3 of （Gohberg et al., 1997）.

The reason is that the problem of computing the coefficients from the roots is not equivalent with the general interpolation problem, as these two problems have completely different condition. See clause 2.2 below.

For the problem at hand, algorithms P, R and R+ are the only algorithms of $O(n^2)$ time complexity to be considered here.

### 1.3 Related problems and algorithms

#### 1.3.1 Problem B: Symmetric functions

By problem B we denote the computation of symmetric functions of order $m$ and parameters $\omega_1,\cdots,\omega_n$, defined by

(5) $\alpha_m = \sum \prod_{1\le i_1\cdots<i_j\cdots<i_m\le n} \omega_{i_j}$

As

(6) $\alpha_m = (-1)^{n-m} a_{n-m}$,

they are obtained from algorithm P, R and R+ just by a change of sign and order. Thus problem A and B are equivalent and have the same solutions.

#### 1.3.2 Problem C: Coefficients of reduced polynomials

By problem C we denote the computation of the coefficients of $P_{\setminus k}(\omega)$

(7) $P_{\setminus k}(\omega) = \prod_{j\ne k} \omega - \omega_j = \frac{P(\omega)}{\omega - \omega_k} = -\frac{\partial P(\omega)}{\partial \omega_k} = \sum a_{m\setminus k}\omega^m$

from the roots and the coefficients of $P(\omega)$. This is based on recursion

(8) $a_{n-1\setminus k} = 1, a_{m-1\setminus k} = a_m + a_{m\setminus k}\omega_k$,

which follows from $P(\omega) = P_{\setminus k}(\omega)(\omega - \omega_k)$ and underlies formula (3). $P_{\setminus k}(\omega)$ and $a_{m\setminus k}$ denote the reduced polynomials and symmetric functions, respectively obtained from omitting the $k$-th root.

Algorithm S2: Solution of problem C using recursion

Input: $k, a_0, \cdots, a_n, \omega_1, \cdots, \omega_n$

Ouput: $a_{0\setminus k}, \cdots, a_{n-1\setminus k}$

Computation:

    a)  $a_{n-1\setminus k} = 1$
    b)  for $i = n-2, \cdots, 0$ do

$$a_{i\setminus k} = a_{i+1} + \omega_k * a_{i+1\setminus k}$$

    end

Algorithm S2 takes $\binom{n}{2} + O(n)$ complex additions and multiplications and $O(n)$ complex storage. The leading term amounts to $4n^2$ floating point operations, reducing to $n^2$ floating point operations in case of real roots.

### 1.3.3 Problem D: Inversion of Vandermondians

By problem D we denote the computation of the inverse of Vandermondian $V = \left(\omega_i^j\right)$ from the parameters. From the extensive literature on the inversion of Vandermondians we reference （Gohberg et al., 1997）, （Parker, 1964） （Traub, 1966）. （Gohberg et al., 1997） evaluated various algorithms including general-purpose and the Björck-Pereyra algorithm, concluding that the Parker-Traub algorithm was the best choice, combining time complexity $O(n^2)$ with the highest accuracy. Their numerical tests were based on real-valued parameters, whereas we will focus on complex parameters.

The inverse Vandermondian is given by

(9)   $V^{-1} = \left(\dfrac{a_{i\setminus j}}{P_{\setminus j}(\omega_j)}\right),$

where $P_{\setminus k}(\omega) = \prod_{j \neq k} \omega - \omega_j$. Evidently, $\sum_{j=0}^{n-1} \dfrac{a_{j\setminus k}}{P_{\setminus k}(\omega_k)} \omega_i^j = \delta_{i-k}$.

The Parker-Traub algorithm solves problem C as follows:

    a)  Obtain $a_m$ $(0 \leq j \leq n)$ using algorithm R.
    b)  For $k = 1, \cdots, n$ execute algorithm S2, obtaining $a_{i\setminus k}$ $(0 \leq i < n, 1 \leq k \leq n)$.
    c)  Obtain the matrix elements evaluating formula (9).

Again Leja ordering, i.e. replacing algorithm R by R+, greatly improves the accuracy. We denote by algorithm PT+ the combination of algorithm L followed by the Parker-Traub procedure.

Algorithm PT+ has time complexity $O(n^2)$ and storage complexity $O(n)$. The leading term of the time complexity amounts to $\frac{39}{2}n^2$ floating point operations for complex roots, reducing to $6n^2$ floating point operations in case of real roots.

As step a) solves problem A, algorithm R can be replaced by algorithm P. This algorithm for problem D is denoted by PP. Our numerical tests show that Vandermondians of dimension 2'000 and more with roots on the unit circle can be accurately inverted based by algorithm PP. The leading term of the time complexity amounts to $20n^2$ floating point operations, reducing to $7n^2$ floating point operations in case of real roots.

### 1.3.4 Problem E, F, G and H

The letter E is assigned to the problem of solving a linear system with Vandermondian coefficients, given its roots. This corresponds to the solution of problem D followed by matrix multiplication, and should be distinguished from what is commonly understood by solving a linear system with Vandermondian coefficients, as we will see below.

We denote by problem F the evaluation of a polynomial at a point, given its roots. This is simply achieved by $n$ additions and multiplications evaluating the product representation of the polynomial. We will use it, however, to assess the performance of algorithm P, R and R+ in cases where accurate values of the coefficients are not available. There we will measure the difference between the values obtained from the product representation with those obtained from the coefficient representation.

Problem F should be distinguished from the general interpolation problem, denoted by G, where a polynomial is given by general data and such a simple algorithm does not exist.

Problem H is the evaluation of polynomials $P_{\setminus k}(\omega)$ at the roots, where only $P_{\setminus k}(\omega_k) = P'(\omega_k)$ is non-zero. This problem is a special case of Problem F. It will be used to assess the performance of algorithms PT+ and PP in cases where accurate values of the coefficients are not available.

### 1.3.5 Problem I: Coefficients of interpolation polynomials

Problem I is the computation of the coefficients of an interpolation polynomial given by a set of points $(x_i, y_i)$. Lagrange interpolation obtains the Barycentric weights in $O(n^2)$ arithmetic operations, followed by $4n$ arithmetic operations per interpolation value. Possession of the coefficients of the interpolation polynomial halves the cost per interpolation value. Though the cost to obtain the coefficients is higher it may sometimes be worth the investment.

Problem I can be solved from a linear system with Vandermondian coefficients, combining algorithms for problems D and E.

> <u>Algorithm DE: Solution of problem I based on a linear system with Vandermondian coefficients</u>
>
> Input: Points $(x_i, y_i)$ $(0 \leq i \leq n)$
>
> Output: Coefficients of the interpolation polynomial through these points.
>
> Computation:
>
> a) Find the inverse of the Vandermondian with parameters $x_i$, using either algorithm PP or PT
> b) Compute and output $(a_0, \cdots, a_n) = V^{-1}(y_0, \cdots, y_n)$

The computational complexity of algorithm DE equals the sum of the complexities of algorithm PP and PT+, respectively, plus the product of a matrix with a vector. The leading term of the complexity when using algorithm PT amounts to $\frac{55}{2}n^2$ floating point operations for complex data, reducing to $8n^2$ floating point operations in case of real roots. With algorithm PP this becomes $28n^2$ and $9n^2$, respectively. The storage complexity is $O(n^2)$

We propose here an alternative algorithm combining problem G and A.

> <u>Algorithm GA: Solution of problem I based on interpolation and algorithm P</u>
>
> Input: Points $(x_i, y_i)$ $(0 \leq i \leq n)$
>
> Output: Coefficients of the interpolation polynomial through these points.
>
> Computation:
>
> a) Find the Barycentric weights $w_j = y_j / P_{\setminus j}(x_j)$ of Lagrange interpolation, where $P(x) = \prod x - x_i$.
> b) Find the interpolation values $p_j = P(\omega^j) \sum_k \frac{w_k}{\omega^j - x_k}$ at the $n-th$ unit roots.
> c) $a_0, \cdots, a_n 0_{n+1}, \cdots, 0_{N-1} = FFT(p_0, \cdots, p_{N-1})$

d) Output $a_0, \cdots, a_n$

In terms of numerical performance algorithm GA is superior over the conventional approach (see clause 6).

The computational complexity of algorithm GA equals the sum of the complexities of calculating the Barycentric weights plus $n$ complex interpolation values plus the FFT. For complex data algorithm GA takes $3n^2$ complex additions and multiplications. The leading term of the time complexity amounts to $24n^2$ floating point operations for complex roots, reducing to $\frac{19}{2}n^2$ floating point operations in case of real roots. In the latter case the symmetry $p_{n-j} = \bar{p}_j$ halves the number of interpolations, and the division of a real by a complex number is achievable by 6 floating point operations. The storage complexity is $O(n)$.

Algorithm DE is slightly slower on complex data and slightly faster on real data, but its storage complexity is $O(n^2)$, while algorithm GA requires only $O(n)$ storage.

## 2 Problems and conditions

### 2.1 Evaluation methodology

The condition $\kappa_f(x)$ of a function $f(x)$ at point $x$ where $f(x) \neq 0$ equals $\lim_{\xi \to x} \left| \frac{f'(\xi)x}{f(\xi)} \right|$. It is however of no use at the zeros, due to the following:

<ins>Proposition 1</ins>

If $x \neq 0$ is a zero of $f$, and if $f$ is analytic and non-vanishing in a neighborhood of $x$, then $\lim_{\xi \to x} \left| \frac{f'(\xi)x}{f(\xi)} \right|$ is either undefined or infinite.

Proof: If $f(x) = f'(x) = 0$, then the De l'Hôpital rule yields $\frac{f'(x)}{f(x)} = \frac{f''(x)}{f'(x)}$, so the r.h.s. has again a zero denominator. If $f$ is a polynomial we eventually obtain a non-zero numerator divided by a zero denominator.

As polynomials may have some zero coefficients, and the Vandermondian some zero matrix elements, we will consider the condition of vector- and matrix-valued functions guaranteed to yield a non-zero result, rather than the condition of single coefficients and matrix elements, using the following

<ins>Definition 1:</ins> The condition of function $f(d, x)$ with respect to $x$ at $(d, x)$ equals

$$(10)\; \kappa_{f,p}(d;x) = \frac{\left\lVert \left( \frac{\partial f_i(d,x)}{\partial x_j} \right) \right\rVert_p \lVert x \rVert_p}{\lVert f(d,x) \rVert_p},$$

where $\left\lVert \left( \frac{\partial f_i(d,x)}{\partial x_j} \right) \right\rVert_p = \sup_{\lVert \delta \rVert_p = 1} \left\lVert \sum \frac{\partial f_i(d,x)}{\partial x_j} \delta_j \right\rVert_p$, and $\lVert y \rVert_p = (\sum |y_i|^p)^{\frac{1}{p}}$ for vector $y$.

The condition of composite functions satisfies

<ins>Proposition 2: Condition of composite functions</ins>

$\kappa_{f \circ g, p}(x) \leq \kappa_{f,p}(g(x)) \kappa_{g,p}(x)$

Proof:

$$\kappa_{f \circ g, p}(x) = \frac{\sup_{\lVert \delta \rVert_p = 1} \left\lVert \sum \frac{\partial f_i(g(x))}{\partial x_k} \delta_k \right\rVert_p \lVert x \rVert_p}{\lVert f(g(x)) \rVert_p} = \frac{\sup_{\lVert \delta \rVert_p = 1} \left\lVert \sum \frac{\partial f_i(g(x))}{\partial g_j} \frac{\partial g_j(x)}{\partial x_k} \delta_k \right\rVert_p \lVert x \rVert_p}{\lVert f(g(x)) \rVert_p} \leq$$

$$\leq \frac{\sup_{\lVert \delta' \rVert_p = 1} \left\lVert \sum \frac{\partial f_i(g(x))}{\partial g_j} \delta'_j \right\rVert_p \sup_{\lVert \delta \rVert_p = 1} \left\lVert \frac{\partial g_j(x)}{\partial x_k} \delta_k \right\rVert_p \lVert x \rVert_p}{\lVert f(g(x)) \rVert_p} \frac{\lVert g(x) \rVert_p}{\lVert g(x) \rVert_p} = \kappa_{f,p}(g(x)) \kappa_{g,p}(x) \text{ . q.e.d.}$$

In our numerical tests we will determine the relative errors and mean square sample errors

$$(11)\; \varepsilon_2(x) = \frac{\lVert \varphi(x) - f(x) \rVert_2 \lVert x \rVert_2}{\lVert f(x) \rVert_2}, \; \varepsilon o = \sqrt{\frac{\sum_d \varepsilon_2^2(d)}{s}},$$

where $\varphi$ denotes the function value obtained by the algorithm at hand, $f$ the best known value, $d$ the data, $s$ the sample size, $x$ the parameters.

Our implementation was made in floating point arithmetic with a 52-bit mantissa and a 11-bit exponent. The largest floating point number equals approximately 1.8E+308.

### 2.2 Condition by problem

Here the conditions of problems introduced beforehand are determined. Only problem A and B have the same condition, while all other problems have pairwise different conditions, showing that they are inequivalent. Investigating the asymptotic properties of the conditions under scaling we find that problem

E is uniquely ill-conditioned, much worse than all other problems, and thus neither equivalent with problem A nor with problem F or G.

### 2.2.1 Problem A and B: Polynomial coefficients and symmetric functions

Formula (8) implies $\frac{\partial}{\partial \omega_k} a_m = a_{m \setminus k}$ ($0 \leq m < n$) and thus

$$(12) \kappa_{A,\infty}(\omega_1, \cdots, \omega_n) = \frac{\max_k |\omega_k|}{\max_m |a_m|} \max_m \sum_k |a_{m \setminus k}|$$

Problem B is equivalent.

### 2.2.2 Problem C: Reduced polynomials and symmetric functions

The numerical properties of algorithm S2 are related to those of algorithm R, because the latter can be understood as an $n$-fold composition of the former with opposite signs of the roots.

As D is the combination of A and C, proposition 2 yields $\kappa_{C,\infty} \geq \kappa_{D,\infty}/\kappa_{A,\infty}$.

### 2.2.3 Problem D: Inversion of Vandermondians

By formulae (7), (8) and (9) the matrix elements have the partial derivatives

$$(13) \frac{\partial}{\partial \omega_k}\left(\frac{a_{i \setminus j}}{P_{\setminus j}(\omega_j)}\right) = \frac{1}{P_{\setminus j}(\omega_j)} \begin{cases} -a_{i \setminus jk} + \frac{a_{i \setminus j}}{\omega_j - \omega_k} = \frac{a_{i \setminus k}}{\omega_j - \omega_k} & (j \neq k) \\ a_{i \setminus k} \sum_{l \neq k} \frac{1}{\omega_k - \omega_l} & (j = k) \end{cases}$$

yielding

$$(14) \kappa_{D,\infty}(\omega_1, \cdots, \omega_n) = \frac{\max_k |\omega_k|}{\max_{i,j} \left|\frac{a_{i \setminus j}}{P_{\setminus j}(\omega_j)}\right|} \max_{i,j} \sum_k \left|\frac{\partial}{\partial \omega_k}\left(\frac{a_{i \setminus j}}{P_{\setminus j}(\omega_j)}\right)\right|$$

### 2.2.4 Problem E: Solution of Vandermondian systems

In problem E we consider the function $(d, \omega) = (V(\omega_1, \cdots, \omega_n))^{-1} \cdot d$. The condition with respect to $\omega$ at $(d, \omega)$ equals

$$(15) \kappa_{E,\infty}(d; \omega) = \frac{\max_k |\omega_k|}{\max_i \left|\sum_j \frac{a_{i \setminus j}}{P_{\setminus j}(\omega_j)} d_j\right|} \max_i \sum_k \left|\sum_j \frac{\partial}{\partial \omega_k}\left(\frac{a_{i \setminus j}}{P_{\setminus j}(\omega_j)}\right) d_j\right|.$$

In contrast, the condition $\kappa_{E,\infty}(\omega; d)$ with respect to $d$ equals $\frac{\sup_{\|\delta\|_p=1} \|V^{-1} \cdot \delta\|_p \|d\|_p}{\|V^{-1} \cdot d\|_p}$ which yields $\sup_d \kappa_{E,\infty}(\omega; d) = \lambda_{max}/\lambda_{min}$, the well-known matrix condition of $V$.

### 2.2.5 Problem F: Polynomial evaluation

We consider the . The condition is obtained from $\frac{\partial}{\partial \omega_k} P(d) = -P_{\setminus k}(d) = -\frac{P(d)}{d - \omega_k}$ ($d \neq \omega_k$) as

$$(16) \kappa_{F,\infty}(d; \omega) = \frac{\max_k |\omega_k|}{\max_i |P(d_i)|} \max_j \sum_k \left|\frac{P(d_j)}{d_j - \omega_k}\right|.$$

### 2.2.6 Problem G: Interpolation with general data

The function considered here is $P(d)$, where $P$ is given by $y_0 = P(x_0), \cdots, y_n = P(x_n)$.

In terms of the Lagrange polynomials

$$(17) l_j(d) = \frac{P_{\setminus j}(d)}{P_{\setminus j}(x_j)}$$

the interpolated value equals

$$(18) P(d) = \sum y_j l_j(d).$$

We find

$$(19) \frac{\partial}{\partial y_i} P(d) = l_i(d)$$
$$(20) \frac{\partial}{\partial x_i} P(d) = \sum_{j \neq i} \frac{-y_j l_j(d)}{d - x_i} + \frac{y_j l_j(d) - y_i l_i(d)}{x_j - x_i}$$

The condition $\kappa_{G,\infty}(d; x, y)$ follows from formula (10).

### 2.2.7 Problem H: Polynomial evaluation at the roots

We consider the matrix-valued function $f(\omega) = \left(P_{\setminus k}(\omega_j)\right)$. As the data equal the roots, the condition includes the derivatives with respect to the data. $P_{\setminus k}(\omega_j) = 0$ ($j \neq k$) implies $\frac{\partial}{\partial \omega_l} P_{\setminus k}(\omega_j) = -P_{\setminus kl}(\omega_j) = 0$ ($l \neq j$), so $l = j$ and

$$(21) \frac{\partial}{\partial \omega_j}\left(P_{\setminus k}(\omega_j)\right) = \frac{d}{d\omega}\left(P_{\setminus k}(\omega)\right)\Big|_{\omega=\omega_j} \quad (1 \leq j \leq n).$$

The condition becomes

$$(22) \kappa_{H,\infty}(\omega) = \frac{\max|\omega_k|}{\max_k|P_{\setminus k}(\omega_k)|} \max_{1 \leq j \leq n} \frac{d}{d\omega}\left(P_{\setminus k}(\omega)\right)\Big|_{\omega=\omega_j}.$$

### 2.2.8 Problem I: Coefficients of interpolation polynomials

As problem I is a special case of problem E and the combination of problem G and A, we find the upper bound $\kappa_{I,\infty} \leq \min(\kappa_{A,\infty}\kappa_{G,\infty}, \kappa_{E,\infty})$.

## 2.3 Discussion of conditions

### 2.3.1 Scaling sensitivity

Noticing the trivial facts

> Proposition 3
>
> When applying the scaling transform $x \mapsto \sigma x$ to the data and the roots, then $P(d) \mapsto \sigma^n P(d)$ and $a_m \mapsto \sigma^{n-m} a_m$, it is of interest to discuss the sensitivity of the conditions to scale change.

As $\sigma \to \infty$, in problem A $\max_m|a_m \sigma^{n-m}| = |a_0 \sigma^n|(1 + o(1))$, $\max_m|a_{m \setminus k} \sigma^{n-1-m}| = |a_{0 \setminus k} \sigma^{n-1}|(1 + o(1))$, $\max_k|\omega_k \sigma| = O(\sigma)$, thus the condition $\kappa_A \in O(1)$. Problem B is equivalent and has the same condition. For problem D and E similar arguments lead to the same result. The condition of problem F, G and H is scaling-invariant, because they are defined in terms of polynomials, thus as well $O(1)$. Problem C has an $O(1)$ lower bound, problem I an $O(1)$ upper bound.

As $\sigma \to 0$, in problem A $\max_{m}|a_m \sigma^{n-m}| = |a_n \sigma^0| = 1 = \max_{m}|a_{m\setminus k} \sigma^{n-1-m}|$, but still $\max_{k}|\omega_k \sigma| = O(\sigma)$, so $\kappa_A \in O(\sigma)$. For problem D and E similar arguments however yield $O(1)$. The condition of problems F, G and H is scaling-invariant and thus $O(1)$. Problem C has condition $\kappa_C \in \Omega(\sigma^{-1})$, problem I $\kappa_I \in O(\sigma)$.

### 2.3.2 Problems A, F and H

Proposition 4:

Let $\max_{k}|\omega_k| \stackrel{\text{def}}{=} \rho = 1$. Then the following holds:

a) For any polynomial $1 \le \max_{m}|a_m| \le \binom{n}{n/2}\rho^{n/2}$, $1 \le \max_{m}|a_{m\setminus k}| \le \binom{n-1}{n/2}\rho^{n/2}$. the largest coefficients if all $\omega_k = -\rho$, thus $P(\omega) = (\omega + \rho)^n$. Proof: Directly from formula (3) and the Binomial theorem.

b) The condition of problem A satisfies

(23) $m'\rho \le \frac{\max_{m} m|a_m|}{|a_{m'}|}\rho \le \kappa_A$, where $|a_{m'}| = \max_{m}|a_m|$.

Proof: The l.h.s. of (23) follows from $m'|a_{m'}| \le \max_{m} m|a_m|$, as $|a_{m'}|$ is maximal. Noticing

(24) $P'(\omega) = P(\omega) \sum_k \frac{1}{\omega - \omega_k} = -\sum_k \frac{\partial P(\omega)}{\partial \omega_k}$

and thus

(25) $ma_m = \sum_k a_{m-1\setminus k}$,

the r.h.s. of (23) follows from (25) and the triangle inequality.

c) If $a_1 \ne 0$, then, $\kappa_A = \rho \sum \omega_k^{-1} = const$ for all scaling factors $\sigma > \sigma_0$ for some $\sigma_0$.
Proof: Notice $|a_1| = \max_{m} \sum_k |a_{m-1\setminus k}|$ for all $\sigma > \sigma_1$ for some $\sigma_1$, as $a_m \in O(\sigma^{n-m})$. By formula (8) $a_0 = \sum a_{0\setminus k}\omega_k$ which with formula (25) yields $a_1 = \sum a_{0\setminus k} = a_0 \sum \omega_k^{-1}$. As the absolute value of the sum on the r.h.s. becomes $< 1$ for some $\sigma_0 > \sigma_1$, $|a_0| \ge |a_1|$ for $\sigma > \sigma_0$. q.e.d.

d) If $P(\omega) = (\omega + \rho)^n$ then $\left|m' - \frac{n}{1+\rho}\right| < 1$ and $\kappa_A \approx \frac{n\rho}{\rho+1}$, which increases monotonically from zero to $n$ as $\rho$ grows, but independent from the scale $\sigma$.
Proof: From b) and c) above.

e) If $P(\omega) = \omega^n - \rho^n$ then $\max_{m}|a_m| = \rho^n$, $\max_{m}|a_{m\setminus k}| = \rho^{n-1}$, thus $\kappa_A = n$ independent from $\rho$.

f) The condition of problem F satisfies

(26) $\kappa_F \le \rho \max_{j} \sum_{j \ne k} \left|\frac{1}{d_j - \omega_k}\right|$ independent from the scale $\sigma$.

Proof: For F immediately from $\kappa_F \le \rho \frac{\max_{i}|P(d_i)|}{\max_{i}|P(d_i)|} \max_{j \ne k} \sum_k \left|\frac{1}{d_j - \omega_k}\right|$.

a) The condition of problem H satisfies

(27) $\kappa_H \le \rho \max_{k} \sum_{l \ne k} \left|\frac{1}{\omega_k - \omega_l}\right|$ independent from the scale $\sigma$.

Proof: $\left|\frac{d}{d\omega}(P_{\setminus k}(\omega))\right|_{\omega = \omega_j} \le \max_{j}|P_{\setminus k}(\omega_j)| \max_{k} \sum_{l \ne k} \left|\frac{1}{\omega_k - \omega_l}\right|$.

### 2.3.3 Scale optimization

We will see that the numerical performance of algorithms P, R and R+, and thus also PP and PT, depends much on the diameter of the set of roots in the complex plane. Therefore we will perform tests with various diameters. Throughout the paper the letter $\rho$ will denote the radius.

The dependency on $\rho$ indicates that better results may be obtained by scaling the roots and data into the optimal range by $\sigma$ and re-scaling the results back by $\sigma^{-1}$, according to proposition 2 above. We therefore define for any of these algorithms

<u>Algorithm $X_s$: Algorithm X with scaling and re-scaling</u>

- a) Choose an appropriate scaling factor, $\sigma$.
- b) Apply the scaling transform $x \mapsto \sigma x$ to the data and the roots.
- c) Execute algorithm X.
- d) Re-scale the results by $a_m \mapsto \sigma^{m-n} a_m$ and $P \mapsto \sigma^{-n} P$.

We will find that for $\rho < 1$ algorithm P (and thus PP) is improved by choosing $\sigma = \rho^{-1}$, whereas scaling does not pay for algorithms R, R+ and thus PT.

## 3 Performance on problem A: Polynomial coefficients

### 3.1 Roots

In this clause we determine the performance of algorithms P, R and R+ on instances of problem A, using polynomials whose coefficients we can exactly compute. The starting points are

(28) $P_{1,p,q}(x) = (x^p - 1)^q = \sum_{0 \le k \le q} (-x^p)^{q-k} \binom{q}{k}$ of degree $n = pq$.

(29) $P_{2,p,q}(x) = \left(\frac{x^{p+1}-1}{x-1}\right)^q = \left(\sum_{0 \le k \le p} x^k\right)^q$ of degree $n = pq$.

Their partial derivatives by the roots are

(30) $\frac{\partial}{\partial \omega_i} P_{1,p,q}(x) = q P_{1,p,q-1}(x) \sum_{0 \le k < p} \omega_i^{p-1-k} x^k$.

(31) $\frac{\partial}{\partial \omega_i} P_{2,p,q}(x) = q P_{1,p,q-1}(x) \sum_{0 \le k < p} \frac{1-\omega_i^{p-k}}{1-\omega_i} x^k$.

We will apply scaling to the roots of these polynomials and obtain accurate values for the coefficients and the condition using the above formulae and proposition 2.

The roots are unit roots of order $p$ and $p + 1$, respectively, multiplied by the scaling factor $\sigma$. Roots randomly distributed on and inside circles in the complex plane are subject of clause 4 and 5.

### 3.2 Leja ordering

The following table shows the average errors $\varepsilon_2$ of algorithms P and R with and without Leja ordering on the roots of $P_{1,n,1}$ with $\rho = 1$. Leja ordering is most effective on algorithm R, but ineffective on algorithm P. The same is the case on other choices of the roots.

| | Effect of Leja ordering | | | |
|---|---|---|---|---|
| | P+ | P | R+ | R |
| n | Leja on | Leja off | Leja on | Leja off |
| 10 | 1.86E-15 | 1.76E-15 | 1.14E-15 | 1.03E-14 |
| 20 | 3.89E-15 | 3.93E-15 | 3.29E-15 | 1.51E-12 |
| 30 | 8.93E-15 | 8.92E-15 | 3.06E-14 | 4.71E-10 |
| 40 | 7.51E-15 | 7.55E-15 | 7.28E-14 | 3.31E-07 |
| 50 | 1.43E-14 | 1.43E-14 | 6.79E-14 | 4.83E-05 |
| 60 | 1.55E-14 | 1.55E-14 | 1.28E-12 | 1.56E-02 |
| 70 | 1.26E-14 | 1.26E-14 | 3.33E-11 | 4.50E+00 |

We will therefore drop algorithm R without Leja ordering, omit Leja ordering on algorithm P and subsequently compare algorithm P (without Leja ordering) and algorithm R+ (with Leja ordering).

### 3.3 Effect of the spread

The next two tables show the average errors $\varepsilon_2$ of algorithms P and R for various values of the radius $\rho$ on the roots of $P_{1,n,1}$. NaN ("not a number") means that $\rho^{2n}$ exceeds the domain of floating point numbers. This happens when evaluating $\|x\|_2$ during the calculation of $\varepsilon_2$ by formula (11), but not during the execution of either algorithm!

| | Algorithm P for various values of $\rho$ | | | | | | |
|---|---|---|---|---|---|---|---|
| n | 1 | 1.2 | 1.4 | 1.5 | 0.9 | 0.5 | 0.1 |
| 10 | 1.86E-15 | 1.34E-15 | 1.23E-15 | 1.28E-15 | 2.03E-15 | 2.23E-15 | 2.07E-15 |
| 110 | 2.82E-14 | 1.68E-14 | 1.66E-14 | 1.58E-14 | 2.61E-14 | 2.59E-14 | 2.61E-14 |
| 210 | 5.15E-14 | 2.86E-14 | 2.86E-14 | 2.85E-14 | 4.48E-14 | 4.53E-14 | 4.57E-14 |
| 310 | 7.37E-14 | 3.83E-14 | 3.79E-14 | 3.76E-14 | 6.38E-14 | 6.47E-14 | 6.55E-14 |
| 410 | 1.02E-13 | 4.69E-14 | 4.57E-14 | 4.49E-14 | 9.35E-14 | 9.29E-14 | 9.40E-14 |
| 510 | 1.31E-13 | 7.55E-14 | 7.24E-14 | 7.16E-14 | 1.26E-13 | 1.26E-13 | 1.27E-13 |
| 610 | 1.41E-13 | 8.12E-14 | 7.88E-14 | 7.79E-14 | 1.12E-13 | 1.14E-13 | 1.16E-13 |
| 710 | 1.70E-13 | 9.69E-14 | 9.64E-14 | 9.43E-14 | 1.41E-13 | 1.42E-13 | 1.44E-13 |
| 810 | 2.07E-13 | 1.07E-13 | 1.06E-13 | 1.04E-13 | 1.84E-13 | 1.85E-13 | 1.87E-13 |
| 910 | 2.35E-13 | 1.19E-13 | 1.15E-13 | NaN | 2.11E-13 | 2.13E-13 | 2.16E-13 |
| 1010 | 2.67E-13 | 1.29E-13 | 1.25E-13 | | 2.43E-13 | 2.45E-13 | 2.47E-13 |
| 1110 | 2.55E-13 | 1.45E-13 | Nan | | 1.90E-13 | 1.93E-13 | 1.97E-13 |
| 1210 | 2.81E-13 | 1.59E-13 | | | 2.18E-13 | 2.22E-13 | 2.26E-13 |
| 1310 | 3.25E-13 | 1.84E-13 | | | 2.77E-13 | 2.80E-13 | 2.84E-13 |
| 1410 | 3.62E-13 | 1.87E-13 | | | 3.00E-13 | 3.03E-13 | 3.07E-13 |
| 1510 | 3.88E-13 | 1.91E-13 | | | 3.37E-13 | 3.40E-13 | 3.44E-13 |
| 1610 | 4.15E-13 | 2.12E-13 | | | 3.60E-13 | 3.63E-13 | 3.68E-13 |
| 1710 | 4.49E-13 | 2.16E-13 | | | 4.00E-13 | 4.04E-13 | 4.08E-13 |
| 1810 | 4.72E-13 | 2.38E-13 | | | 4.24E-13 | 4.28E-13 | 4.33E-13 |
| 1910 | 5.04E-13 | 2.55E-13 | | | 4.58E-13 | 4.62E-13 | 4.67E-13 |
| 2010 | 5.20E-13 | NaN | | | 4.82E-13 | 4.87E-13 | 4.92E-13 |

Algorithm P is fairly accurate for a range of radii and problem sizes. The accuracy does not strongly depend from the radius, but is highest at $\rho = 1.4$.

| | Algorithm R+ for various values of $\rho$ | | | | | | |
|---|---|---|---|---|---|---|---|
| n | 1 | 1.2 | 1.4 | 1.5 | 0.9 | 0.5 | 0.1 |
| 10 | 1.14E-15 | 1.25E-15 | 1.29E-15 | 1.24E-15 | 6.23E-16 | 1.82E-16 | 1.99E-17 |
| 110 | 2.75E-09 | 7.17E-14 | 1.70E-14 | 1.63E-14 | 5.36E-08 | 6.75E-03 | 3.73E-14 |
| 210 | 1.26E-03 | 2.19E-11 | 2.90E-14 | 2.92E-14 | >1 | >1 | 2.34E-11 |
| 310 | >1 | 8.74E-11 | 3.86E-14 | 3.68E-14 | | | 9.23E-09 |
| 410 | | 1.05E-08 | 1.16E-13 | 5.51E-14 | | | 6.04E-06 |
| 510 | | 4.31E-02 | 4.70E-12 | 7.22E-14 | | | 2.83E-03 |
| 610 | | 3.89E-04 | 7.84E-13 | 2.98E-13 | | | 9.13E-01 |
| 710 | | >1 | 9.21E-12 | 4.30E-13 | | | >1 |
| 810 | | | 4.49E-12 | 1.47E-13 | | | |
| 910 | | | 2.90E-10 | >1 | | | |
| 1010 | | | 3.44E-09 | | | | |
| 1110 | | | NaN | | | | |

Algorithm R+ is more accurate than algorithm P for $n = 10$, but becomes rapidly less accurate as the problem size increases. It achieves the best accuracy and also the largest problem size at $\rho = 1.4$.

### 3.4 Effect of scaling

The previous tables indicate that scaling, as specified in clause 2.3.12.3.3 above, might improve the results. To investigate this we compare algorithm P (without scaling) and algorithm $P_S$ (with scaling), choosing $\sigma = 1/\rho$, for some values of $\rho < 1$.

| | Effect of scaling on algorithm P | | | | | |
|---|---|---|---|---|---|---|
| | P | $P_s$ | P | $P_s$ | P | $P_s$ |
| n | 0.9 | 0.9 | 0.5 | 0.5 | 0.1 | 0.1 |
| 10 | 2.03E-15 | 1.36E-15 | 2.23E-15 | 2.39E-16 | 2.07E-15 | 5.61E-17 |
| 110 | 2.61E-14 | 5.94E-15 | 2.59E-14 | 1.40E-15 | 2.61E-14 | 2.57E-16 |
| 210 | 4.48E-14 | 5.21E-15 | 4.53E-14 | 1.19E-15 | 4.57E-14 | 1.75E-16 |
| 310 | 6.38E-14 | 5.34E-15 | 6.47E-14 | 1.83E-15 | 6.55E-14 | 3.38E-16 |
| 410 | 9.35E-14 | 1.73E-14 | 9.29E-14 | 2.73E-15 | 9.40E-14 | 5.57E-16 |
| 510 | 1.26E-13 | 1.87E-14 | 1.26E-13 | 3.64E-15 | 1.27E-13 | 6.30E-16 |
| 610 | 1.12E-13 | 9.71E-15 | 1.14E-13 | 1.79E-15 | 1.16E-13 | 4.65E-16 |
| 710 | 1.41E-13 | 1.44E-14 | 1.42E-13 | 2.52E-15 | 1.44E-13 | 4.84E-16 |
| 810 | 1.84E-13 | 1.16E-14 | 1.85E-13 | 3.60E-15 | 1.87E-13 | 7.57E-16 |
| 910 | 2.11E-13 | 1.08E-14 | 2.13E-13 | 3.35E-15 | 2.16E-13 | 6.66E-16 |
| 1010 | 2.43E-13 | 1.45E-14 | 2.45E-13 | 5.16E-15 | 2.47E-13 | 1.09E-15 |
| 1110 | 1.90E-13 | 1.24E-14 | 1.93E-13 | 2.28E-15 | 1.97E-13 | 3.84E-16 |
| 1210 | 2.18E-13 | 1.63E-14 | 2.22E-13 | 3.50E-15 | 2.26E-13 | 4.79E-16 |
| 1310 | 2.77E-13 | 1.37E-14 | 2.80E-13 | 5.53E-15 | 2.84E-13 | 1.10E-15 |
| 1410 | 3.00E-13 | 2.48E-14 | 3.03E-13 | 5.32E-15 | 3.07E-13 | 9.45E-16 |
| 1510 | 3.37E-13 | 2.51E-14 | 3.40E-13 | 4.43E-15 | 3.44E-13 | 8.55E-16 |
| 1610 | 3.60E-13 | 1.85E-14 | 3.63E-13 | 4.24E-15 | 3.68E-13 | 8.71E-16 |
| 1710 | 4.00E-13 | 1.46E-14 | 4.04E-13 | 3.01E-15 | 4.08E-13 | 5.76E-16 |
| 1810 | 4.24E-13 | 1.89E-14 | 4.28E-13 | 5.64E-15 | 4.33E-13 | 9.00E-16 |
| 1910 | 4.58E-13 | 1.96E-14 | 4.62E-13 | 6.22E-15 | 4.67E-13 | 1.21E-15 |
| 2010 | 4.82E-13 | 1.75E-14 | 4.87E-13 | 5.84E-15 | 4.92E-13 | 1.29E-15 |

Indeed, for $\rho < 1$ algorithm scaling increases the accuracy of algorithm P by 1-2 orders of magnitude. For algorithm R+ the effect is negligible, because $a_k$ are less precise for small $k$. For $\rho > 1$ scaling does not pay for either algorithm.

Similar but not identical results are obtained for other choices of the roots. Varying $\sigma$ between $0.9/\rho$ and $1.1/\rho$ sometimes achieves the best values for algorithm P, sometimes also for algorithm R+. However, almost optimal is $\sigma = 1/\rho$.

Unless stated otherwise subsequently algorithm P will be used for $\rho \geq 1$, algorithm P$_s$ with $\sigma = 1/\rho$ for $\rho < 1$ and algorithm R+ for all $\rho$.

### 3.5 Performance on unit roots without unity

The following two tables present the relative errors of algorithm P and R+ to determine the coefficients of the polynomial $P_{2,n,1}(x) = \frac{x^n - 1}{x - 1}$ for various values of $\rho$.

| | P: Unit roots without unity | | | | |
|---|---|---|---|---|---|
| n | 1 | 1.2 | 1.4 | 0.5 | 0.1 |
| 10 | 1.92E-15 | 2.25E-15 | 1.56E-15 | 1.01E-15 | 2.13E-16 |
| 110 | 1.64E-14 | 1.15E-14 | 1.04E-14 | 6.62E-15 | 1.33E-15 |
| 210 | 3.41E-14 | 2.38E-14 | 2.24E-14 | 1.40E-14 | 2.70E-15 |
| 310 | 5.78E-14 | 3.79E-14 | 2.98E-14 | 5.69E-14 | 1.11E-14 |
| 410 | 9.08E-14 | 4.03E-14 | 3.69E-14 | 7.53E-14 | 1.51E-14 |
| 510 | 1.20E-13 | 1.37E-13 | 1.38E-13 | 8.46E-14 | 1.72E-14 |
| 610 | 1.47E-13 | 5.69E-14 | 7.53E-14 | 6.46E-14 | 1.27E-14 |
| 710 | 1.48E-13 | 7.11E-14 | 7.65E-14 | 7.37E-14 | 1.50E-14 |
| 810 | 3.34E-13 | 1.14E-13 | 1.11E-13 | 8.72E-14 | 1.72E-14 |
| 910 | 1.51E-13 | 7.77E-14 | 7.34E-14 | 1.02E-13 | 2.05E-14 |
| 1010 | 2.70E-13 | 1.11E-13 | 1.13E-13 | 1.09E-13 | 2.16E-14 |
| 1110 | 2.08E-13 | 1.85E-13 | NaN | 2.02E-13 | 4.03E-14 |
| 1210 | 2.56E-13 | 2.58E-13 | | 2.27E-13 | 4.56E-14 |
| 1310 | 4.38E-13 | 1.17E-13 | | 2.51E-13 | 5.02E-14 |
| 1410 | 3.23E-13 | 1.47E-13 | | 2.72E-13 | 5.36E-14 |
| 1510 | 3.65E-13 | 1.80E-13 | | 2.85E-13 | 5.75E-14 |
| 1610 | 5.94E-13 | 1.87E-13 | | 2.98E-13 | 5.93E-14 |
| 1710 | 4.66E-13 | 1.63E-13 | | 3.10E-13 | 6.20E-14 |
| 1810 | 5.81E-13 | 1.69E-13 | | 3.27E-13 | 6.64E-14 |
| 1910 | 6.34E-13 | 2.78E-13 | | 3.35E-13 | 6.76E-14 |
| 2010 | 4.42E-13 | NaN | | 3.46E-13 | 6.99E-14 |

The picture for both algorithms is similar to that in the previous clause.

| | R+: Unit roots without unity | | | |
|---|---|---|---|---|
| n | 1 | 1.2 | 1.4 | 0.5 |
| 10 | 1.01E-15 | 1.14E-15 | 1.19E-15 | 4.12E-16 |
| 110 | 8.47E-11 | 2.85E-14 | 9.64E-15 | 1.04E-02 |
| 210 | 6.69E-03 | 5.06E-11 | 3.30E-14 | |
| 310 | | 1.04E-10 | 4.44E-14 | |
| 410 | | 1.14E-07 | 8.72E-14 | |
| 510 | | 1.70E-06 | 1.33E-13 | |
| 610 | | 3.04E-04 | 5.17E-13 | |
| 710 | | | 2.19E-11 | |
| 810 | | | 4.90E-12 | |
| 910 | | | 3.03E-11 | |
| 1010 | | | 1.66E-09 | |
| 1110 | | | NaN | |

### 3.6 Effect of $q$

The maximal coefficients of the polynomials $P_{1,p,q}(x)$ and $P_{2,p,q}(x)$ grow exponentially with $q$ while keeping $n = pq$ approximately constant. Algorithm R+ is very sensitive to an increase of $q$, while algorithm P manages much larger values

In the following tables $\rho = 1$, in the left one $q = 6$, in the right one $q = 120$.

| | $P_{2,p,6}$ | |
|---|---|---|
| n | P | R+ |
| 12 | 3.37E-15 | 2.18E-15 |
| 114 | 1.58E-14 | 8.82E-02 |
| 210 | … | >1 |
| 2010 | 4.93E-13 | |

| | $P_{2,p,120}$ | |
|---|---|---|
| n | P | R+ |
| 120 | 2.07E-14 | 7.36E-15 |
| 240 | 5.93E-14 | 1.28E-02 |
| 360 | … | >1 |
| 2040 | 4.64E-13 | |

Algorithm P performs excellently even for $P_{1,4,500}$. Algorithm R+ performs increasingly better as $q$ grows, but stays behind algorithm P until $q$ approaches $n$. At $q = n$ both algorithms work excellently, algorithm R+ even better than algorithm P. The tables below compare the results of the two algorithms for $\rho = 1$.

| | $P_{1,1,n}$ | |
|---|---|---|
| n | P | R+ |
| 10 | 8.18E-16 | 0 |
| 110 | 1.04E-14 | 4.93E-16 |
| 210 | 1.98E-14 | 2.16E-16 |
| 310 | 2.84E-14 | 1.89E-16 |
| 410 | 3.96E-14 | 1.11E-15 |
| 510 | 4.99E-14 | 7.90E-16 |
| 610 | NaN | NaN |

| | $P_{2,1,n}$ | |
|---|---|---|
| n | P | R+ |
| 10 | 2.06E-15 | 6.28E-16 |
| 110 | 1.89E-14 | 6.75E-15 |
| 210 | 3.33E-14 | 1.29E-14 |
| 310 | 5.00E-14 | 1.90E-14 |
| 410 | 6.91E-14 | 2.51E-14 |
| 510 | 1.02E-13 | 3.12E-14 |
| 610 | NaN | NaN |

### 3.7 Conclusion

For the solution of problem A algorithm P outperforms algorithm R+ case of unit roots. This has been established for $0 < \rho < 1.4$. The performance of algorithm R+ improves as $\rho \to 1.4$ and as $q \to n$. It is only in the latter case that algorithm R+ beats algorithm P by 1-2 decimals.

## 4 Performance on problem F: Polynomial evaluation

While previously only unit roots have been considered, we turn now to roots randomly sampled from certain distributions. As we have no means to determine the coefficients in a more reliable way than by our algorithms, we will instead of problem A turn to problem F.

Specifically, we choose data $x_0, \cdots, x_n$ and compare $y_i = \prod_k x_i - \omega_k$ with $z_i = \sum_k a_k x_i^k$, where the coefficients are obtained from our algorithms.

In the following $\delta$ denotes a random variable equally distributed in [0,1]. The average errors are obtained according to formula (11) from 100 samples of roots and data for $n < 255$, and from 10 samples otherwise.

While previously the NaN condition arose during evaluation of $\varepsilon_2$ with $\sigma > 1$ due to overflow, now it occurs additionally when all $y_i$ vanish after re-scaling with $\sigma < 1$ due to underflow,

### 4.1 Circle

Roots and data are sampled from the random variables $d_k = \rho e^{i(k+\delta)/n}$ on the complex circle with radius $\rho$. The next two tables show the results of the two algorithms for various $\rho$.

| | P:Circle | | | | |
|---|---|---|---|---|---|
| n | 1 | 1.2 | 1.4 | 0.5 | 0.1 |
| 10 | 1.56E-15 | 1.63E-15 | 1.67E-15 | 1.50E-15 | 1.52E-15 |
| 110 | 1.73E-14 | 1.84E-14 | 1.76E-14 | 1.65E-14 | 1.73E-14 |
| 210 | 3.17E-14 | 3.35E-14 | 3.33E-14 | 3.23E-14 | NaN |
| 310 | 4.28E-14 | 4.79E-14 | 4.32E-14 | 4.33E-14 | |
| 410 | 6.27E-14 | 6.89E-14 | 6.57E-14 | 6.01E-14 | |
| 510 | 8.49E-14 | 8.52E-14 | 9.03E-14 | NaN | |
| 610 | 8.17E-14 | 8.57E-14 | 8.27E-14 | | |
| 710 | 1.00E-13 | 1.04E-13 | 1.00E-13 | | |
| 810 | 1.23E-13 | 1.28E-13 | 1.26E-13 | | |
| 910 | 7.18E-14 | 1.46E-13 | 1.50E-13 | | |
| 1010 | 1.63E-13 | 1.73E-13 | 1.68E-13 | | |
| 1110 | 1.42E-13 | 1.47E-13 | NaN | | |
| … | … | … | | | |
| 1910 | 3.13E-13 | 3.17E-13 | | | |
| 2010 | 3.33E-13 | NaN | | | |

| | R+:Circle | | | | |
|---|---|---|---|---|---|
| n | 1 | 1.2 | 1.4 | 0.5 | 0.1 |
| 10 | 4.67E-16 | 5.67E-16 | 5.53E-16 | 2.63E-15 | 2.34E-15 |
| 110 | 1.04E-08 | 4.68E-12 | 1.09E-13 | >1 | >1 |
| 210 | 5.88E-02 | 1.81E-07 | 1.02E-10 | | |
| 310 | >1 | 8.99E-05 | 1.63E-09 | | |
| 410 | | 4.01E-01 | 7.86E-09 | | |
| 510 | | >1 | 9.21E-06 | | |
| 610 | | | 2.87E-06 | | |
| 710 | | | 3.54E-05 | | |
| 810 | | | 2.32E-01 | | |
| 910 | | | 6.20E-02 | | |
| 1010 | | | >1 | | |

Algorithm P remains exact in the widest possible range, while algorithm R+ becomes unstable beforehand.

## 4.2 Disk

Roots and data are sampled from the random variables $d_k = \delta\rho e^{i(k+\delta)/n}$ in the complex disk with radius $\rho$. The next two tables show the results of the two algorithms for various $\rho$.

| | P:Disk | | | | |
|---|---|---|---|---|---|
| n | 1 | 1.2 | 1.4 | 0.5 | 0.1 |
| 10 | 2.04E-15 | 2.95E-15 | 2.83E-15 | 1.72E-15 | 3.72E-15 |
| 110 | 3.93E-14 | 3.43E-14 | 1.50E-13 | 6.30E-14 | 1.48E-13 |
| 210 | 9.90E-14 | 3.51E-12 | 1.65E-13 | 3.28E-14 | NaN |
| 310 | 2.17E-13 | 2.05E-12 | 3.50E-13 | 6.73E-13 | |
| 410 | 1.99E-13 | 2.10E-12 | 1.10E-09 | 2.46E-13 | |
| 510 | 5.97E-13 | 1.69E-12 | 6.35E-11 | NaN | |
| 610 | 2.19E-13 | 6.65E-12 | 1.20E-11 | | |
| 710 | 1.78E-12 | 1.55E-10 | 1.73E-10 | | |
| 810 | 4.76E-12 | 5.73E-12 | 1.07E-11 | | |
| 910 | 3.32E-12 | 1.19E-09 | 3.36E-10 | | |
| 1010 | 1.39E-11 | 1.27E-09 | 1.82E-09 | | |
| 1110 | 2.15E-12 | 1.78E-12 | NaN | | |
| … | … | … | | | |
| 1910 | 1.91E-11 | 1.73E-10 | | | |
| 2010 | 1.81E-12 | NaN | | | |

| | R+:Disk | | | | |
|---|---|---|---|---|---|
| n | 1 | 1.2 | 1.4 | 0.5 | 0.1 |
| 10 | 5.40E-16 | 5.19E-16 | 5.93E-16 | 1.49E-15 | 6.32E-16 |
| 110 | 3.18E-09 | 5.82E-12 | 2.12E-13 | >1 | 8.65E-04 |
| 210 | 2.08E-02 | 5.01E-05 | 3.87E-09 | | >1 |
| 310 | >1 | >1 | 1.19E-03 | | |
| 410 | | | >1 | | |

Algorithm P remains exact in the widest possible range, while algorithm R+ becomes unstable beforehand.

## 4.3 Annulus

Roots and data are sampled from the random variables $d_k = (1-\Delta\delta)\rho e^{i(k+\delta)/n}$ in the complex disk with radius $\rho$. The next two tables show the results of the two algorithms for various $\rho$ and $\Delta = 0.1$.

| | P:Annulus | | | | |
|---|---|---|---|---|---|
| n | 1 | 1.2 | 1.4 | 0.5 | 0.1 |
| 10 | 1.43E-15 | 1.58E-15 | 1.46E-15 | 1.45E-15 | 1.49E-15 |
| 110 | 7.63E-14 | 7.12E-14 | 6.84E-14 | 6.44E-14 | 8.48E-14 |
| 210 | 1.38E-13 | 9.40E-14 | 1.07E-13 | 1.50E-13 | NaN |
| 310 | 6.71E-13 | 1.50E-13 | 3.16E-13 | 3.38E-13 | |
| 410 | 6.80E-13 | 4.67E-13 | 6.22E-13 | 7.47E-13 | |
| 510 | 1.42E-12 | 5.72E-13 | 1.45E-12 | NaN | |
| 610 | 1.61E-12 | 7.76E-13 | 3.16E-12 | | |
| 710 | 4.50E-12 | 7.51E-12 | 3.12E-11 | | |
| 810 | 1.01E-12 | 1.15E-12 | 2.47E-12 | | |
| 910 | 2.86E-12 | 3.27E-12 | 2.64E-12 | | |
| 1010 | 1.91E-12 | 4.54E-12 | 8.57E-12 | | |
| 1110 | 6.79E-13 | 8.61E-13 | NaN | | |
| … | … | … | | | |
| 1910 | 1.59E-11 | 6.88E-10 | | | |
| 2010 | 2.63E-11 | NaN | | | |

| | R+:Annulus | | | | |
|---|---|---|---|---|---|
| n | 1 | 1.2 | 1.4 | 0.5 | 0.1 |
| 10 | 4.73E-16 | 5.29E-16 | 4.61E-16 | 1.37E-15 | 1.66E-15 |
| 110 | 5.80E-08 | 5.90E-12 | 4.72E-14 | >1 | >1 |
| 210 | >1 | 1.01E-08 | 7.97E-12 | | |
| 310 | | 4.85E-05 | 3.94E-10 | | |
| 410 | | 2.41E-03 | 2.26E-09 | | |
| 510 | | >1 | 1.38E-08 | | |
| 610 | | | 7.16E-07 | | |
| 710 | | | 1.03E-06 | | |
| 810 | | | 2.87E-01 | | |
| 910 | | | 2.13E-04 | | |
| 1010 | | | 3.37E-02 | | |
| 1110 | | | NaN | | |

The picture almost the same as before, R+ behaving somewhat better.

## 4.4 Line

Roots are sampled from the random variables $d_k = \rho(-1 + 2(k+\delta)/n)$ the data from $d_k + \tau$. The next two tables show the results of the two algorithms for various $\rho$ and $\tau = 0$, i.e. data and roots have the same distribution.

| | P:Line | | | | |
|---|---|---|---|---|---|
| n | 1 | 1.2 | 1.4 | 0.5 | 0.1 |
| 10 | 2.27E-15 | 4.00E-15 | 2.26E-14 | 4.63E-12 | 8.90E-04 |
| 20 | 3.20E-14 | 8.33E-13 | 5.97E-11 | 1.91E-06 | >1 |
| 30 | 2.65E-13 | 2.46E-10 | 1.14E-07 | >1 | |
| 40 | 1.62E-12 | 7.71E-08 | 2.53E-05 | | |
| 50 | 6.41E-13 | 5.58E-07 | >1 | | |
| 60 | 1.72E-09 | 9.83E-02 | | | |
| 70 | 1.88E-03 | 3.20E-01 | | | |
| 80 | 2.12E-01 | >1 | | | |
| 90 | 1.60E-03 | | | | |
| 100 | >1 | | | | |

| | R+:Line | | | | |
|---|---|---|---|---|---|
| n | 1 | 1.2 | 1.4 | 0.5 | 0.1 |
| 10 | 5.31E-16 | 1.14E-15 | 5.32E-15 | 3.17E-13 | 1.28E-04 |
| 20 | 3.52E-15 | 5.77E-14 | 4.81E-12 | 7.19E-08 | >1 |
| 30 | 1.56E-14 | 8.49E-12 | 4.85E-08 | 7.52E-03 | |
| 40 | 2.26E-13 | 2.68E-09 | 9.41E-07 | >1 | |
| 50 | 3.08E-14 | 1.93E-05 | 8.13E-02 | | |
| 60 | 2.24E-10 | 1.99E-03 | 3.17E-02 | | |
| 70 | 1.87E-04 | 1.04E-03 | >1 | | |
| 80 | 7.09E-02 | 2.13E-03 | | | |
| 90 | 1.74E-04 | >1 | | | |
| 100 | >1 | | | | |

Unlike in the previous cases, here algorithm R+ yields errors of 1-2 magnitudes smaller than algorithm P. Both algorithms perform much poorer than they do on the disk, the annulus and the circle, both becoming unstable below $n = 100$.

To gain further insight we change the distribution of the data. The next two tables display the results for various $\rho$ and $\tau = 1$, i.e. the data are shifted to the right. Unlike in all other cases we use here algorithm P without scaling for all $\rho$, as the scaling of the shifted data with $\sigma = 1/\rho$ causes overflow.

| | P:Line with shifted data | | | | |
|---|---|---|---|---|---|
| n | 1 | 1.2 | 1.4 | 0.5 | 0.1 |
| 10 | 8.99E-16 | 1.53E-15 | 1.15E-15 | 1.06E-15 | 1.06E-15 |
| 110 | 8.41E-08 | 2.74E-07 | 9.50E-03 | 3.55E-13 | 3.55E-13 |
| 210 | >1 | >1 | >1 | 5.99E-11 | 5.99E-11 |
| 310 | | | | 6.67E-09 | 6.67E-09 |
| 410 | | | | 2.12E-06 | 2.12E-06 |
| 510 | | | | 3.90E-03 | 3.90E-03 |
| 610 | | | | 1.29E-01 | 1.29E-01 |
| 710 | | | | >1 | 1.83E-13 |
| 810 | | | | | 2.46E-13 |
| 910 | | | | | 3.26E-13 |
| 1010 | | | | | 4.30E-13 |
| 1110 | | | | | … |

| | R+:Line with shifted data | | | | |
|---|---|---|---|---|---|
| n | 1 | 1.2 | 1.4 | 0.5 | 0.1 |
| 10 | 4.46E-16 | 3.60E-16 | 3.41E-16 | 4.56E-16 | 3.77E-16 |
| 110 | 8.52E-13 | 6.52E-12 | 4.46E-11 | 6.05E-15 | 2.55E-15 |
| 210 | 1.36E-08 | 4.43E-07 | 1.02E-05 | 2.06E-13 | 4.89E-15 |
| 310 | 3.10E-05 | 9.83E-03 | >1 | 5.63E-11 | 8.87E-15 |
| 410 | >1 | >1 | | 5.05E-03 | 8.04E-15 |
| 510 | | | | >1 | 8.55E-15 |
| 610 | | | | | 1.62E-14 |
| 710 | | | | | 1.40E-14 |
| 810 | | | | | 1.46E-14 |
| 910 | | | | | 1.73E-14 |
| 1010 | | | | | 7.27E-13 |
| 1110 | | | | | … |

Still algorithm R+ performs somewhat better than algorithm P, but both algorithms remain stable at much larger $n$.

There are two differences: Firstly, the values of the polynomials at the data are now much larger than before, so absolute errors have a smaller effect on the relative errors. Secondly, the differences of the data and the roots are much larger than before, so the condition is much smaller.

But why is the performance on the circle, the disk and the annulus so much better? The condition does not suffice as an explanation, because it is proportional to the average point density. The average point density is on the line higher by a factor $\pi$ than on the circle, but the error at $n\pi$ on the disk is still better than at $n$ on the line.

The decisive factor apparently is the average modulus of the data, which causes larger polynomial values at the data, lowering the effects of absolute errors.

## 5 Performance on problem H: Evaluation at the roots

The product of a Vandermondian with its inverse involves the evaluation of the reduced polynomials at the roots: $P_{\backslash k}(\omega_j) = 0 \ (j \neq k)$. In order to compare the suitability of algorithm P and R+ as building block of algorithms PP and PT+ it is sufficient to evaluate $P(\omega_j)$, because $P_{\backslash k}(\omega_j) = \frac{P(\omega_j)}{\omega_k - \omega_j} \ (j \neq k)$, where only the numerator depends from algorithm P and R+. Therefore we consider here problem H, evaluating $P$ at the roots and at the origin.

The average errors are obtained according to formula (11) from 100 samples of roots and data for $n < 255$, and from 10 samples otherwise.

### 5.1 Circle

Roots are sampled from the random variables $d_k = \rho e^{i(k+\delta)/n}$ on the complex circle with radius $\rho$. The next two tables show the results of the two algorithms for various $\rho$.

| | | P:Circle | | | |
|---|---|---|---|---|---|
| n | 1 | 1.2 | 1.4 | 0.5 | 0.1 |
| 10 | 8.47E-15 | 8.72E-15 | 8.79E-15 | 8.82E-15 | 8.70E-15 |
| 110 | 3.59E-13 | 3.66E-13 | 3.60E-13 | 3.64E-13 | 3.63E-13 |
| 210 | 9.12E-13 | 9.19E-13 | 9.01E-13 | 9.04E-13 | NaN |
| 310 | 1.56E-12 | 1.63E-12 | 1.58E-12 | 1.58E-12 | |
| 410 | 2.59E-12 | 2.61E-12 | 2.57E-12 | 2.64E-12 | |
| 510 | 4.08E-12 | 4.00E-12 | 4.02E-12 | NaN | |
| 610 | 4.05E-12 | 4.09E-12 | 4.08E-12 | | |
| 710 | 5.36E-12 | 5.52E-12 | 5.53E-12 | | |
| 810 | 7.15E-12 | 7.21E-12 | 5.51E-12 | | |
| 910 | 9.19E-12 | 8.87E-12 | 8.88E-12 | | |
| 1010 | 1.08E-11 | 1.07E-11 | 1.05E-11 | | |
| 1110 | 3.04E-12 | 9.56E-12 | NaN | | |
| … | … | … | | | |
| 1910 | 2.79E-11 | 2.78E-11 | | | |
| 2010 | 2.98E-11 | NaN | | | |

| | | R+:Circle | | | |
|---|---|---|---|---|---|
| n | 1 | 1.2 | 1.4 | 0.5 | 0.1 |
| 10 | 2.98E-15 | 2.46E-15 | 2.31E-15 | 1.98E-14 | 1.94E-14 |
| 110 | 3.18E-07 | 1.38E-10 | 2.23E-12 | >1 | >1 |
| 210 | >1 | 1.73E-06 | 8.46E-09 | | |
| 310 | | 1.33E-03 | 8.26E-09 | | |
| 410 | | >1 | 5.23E-07 | | |
| 510 | | | 1.60E-04 | | |
| 610 | | | 3.15E-03 | | |
| 710 | | | 3.38E-01 | | |
| 810 | | | 5.55E-01 | | |
| 910 | | | >1 | | |

Algorithm P remains exact in the widest possible range, while algorithm R+ becomes unstable beforehand.

### 5.2 Disk

Roots are sampled from the random variables $d_k = \delta \rho e^{i(k+\delta)/n}$ in the complex disk with radius $\rho$.

| P:Disk | | | | | |
|---|---|---|---|---|---|
| n | 1 | 1.2 | 1.4 | 0.5 | 0.1 |
| 10 | 7.02E-12 | 1.74E-12 | 2.74E-12 | 1.64E-09 | 2.95E-12 |
| 20 | 7.90E-08 | 9.96E-08 | 1.04E-06 | 2.27E-08 | 5.70E-09 |
| 30 | 1.15E-02 | 3.72E-02 | 2.08E-02 | 2.66E-04 | >1 |
| 40 | >1 | >1 | >1 | >1 | |

| R+:Disk | | | | | |
|---|---|---|---|---|---|
| n | 1 | 1.2 | 1.4 | 0.5 | 0.1 |
| 10 | 1.37E-12 | 9.82E-14 | 4.12E-13 | 7.71E-11 | 6.73E-13 |
| 20 | 1.29E-08 | 9.84E-09 | 7.11E-08 | 1.80E-08 | 9.57E-09 |
| 30 | 8.33E-03 | 1.28E-03 | 2.48E-03 | 1.69E-03 | >1 |
| 40 | >1 | >1 | >1 | >1 | |

Both algorithms perform similarly poorly.

Problem H is clearly harder than problem F, because the denominator of $\varepsilon_2$ is much smaller.

### 5.3 Line

Roots are sampled from the random variables $d_k = \delta(-1 + 2(k+\delta)/n)$.

| P:Line | | | | | |
|---|---|---|---|---|---|
| n | 1 | 1.2 | 1.4 | 0.5 | 0.1 |
| 10 | 3.72E-12 | 6.27E-11 | 5.02E-10 | 1.06E-10 | 6.25E-07 |
| 20 | 3.21E-03 | 2.50E-08 | 1.37E-06 | 8.68E-04 | >1 |
| 30 | 4.40E-03 | >1 | 3.91E-01 | >1 | |
| 40 | >1 | | >1 | | |

| R+:Line | | | | | |
|---|---|---|---|---|---|
| n | 1 | 1.2 | 1.4 | 0.5 | 0.1 |
| 10 | 6.39E-13 | 1.10E-11 | 6.40E-11 | 3.39E-12 | 1.21E-13 |
| 20 | 4.27E-06 | 1.27E-08 | 1.02E-06 | 1.78E-06 | 2.26E-10 |
| 30 | 1.31E-04 | >1 | 1.71E-02 | 7.54E-03 | 4.32E-07 |
| 40 | >1 | | >1 | >1 | 5.95E-04 |
| 50 | | | | | >1 |

Both algorithms perform poorly, R+ slightly less than P.

## 6 Performance on problem I: Coefficients of interpolation polynomials

To assess the numerical errors of algorithms DE and GA, we will randomly choose roots and arguments, evaluate the polynomial at the arguments, determine the coefficients with algorithm P from the roots and calculate the errors as the difference them and the coefficients yielded by algorithms DE and GA from the data.

### 6.1 Annulus

Roots and data are sampled from the random variables $d_k = (1 - \Delta\delta)\rho e^{i(k+\delta)/n}$ in the complex disk with radius $\rho$. The next two tables show the results of the two algorithms for various $\rho$ and $\Delta = 0.1$.

| | P: Annulus | | | | |
|---|---|---|---|---|---|
| n | 1 | 1.2 | 1.4 | 0.5 | 0.1 |
| 10 | 5.54E-15 | 3.05E-15 | 7.93E-15 | 4.14E-15 | 4.21E-15 |
| 110 | 6.44E-15 | 1.93E-12 | 1.66E-12 | 2.70E-14 | 2.70E-14 |
| 210 | 7.52E-14 | 4.13E-10 | 4.24E-10 | 4.79E-14 | NaN |
| 310 | 3.46E-14 | 9.62E-08 | 1.41E-07 | 6.82E-14 | |
| 410 | 1.26E-13 | 1.68E-05 | 4.80E-05 | 2.60E-13 | |
| 510 | 2.97E-13 | 7.98E-03 | 5.07E-03 | NaN | |
| 610 | 3.04E-13 | 5.25E-01 | 8.14E-01 | | |
| 710 | 1.78E-12 | >1 | 9.99E-01 | | |
| 810 | 1.46E-11 | | >1 | | |
| 910 | 1.14E-11 | | | | |
| 1010 | 6.42E-11 | | | | |
| 1110 | 7.93E-12 | | | | |
| | … | | | | |
| 1910 | 1.37E-09 | | | | |
| 2010 | 1.67E-10 | | | | |

| | R+:Annulus | | | | |
|---|---|---|---|---|---|
| n | 1 | 1.2 | 1.4 | 0.5 | 0.1 |
| 10 | 8.69E-15 | 1.33E-14 | 1.58E-14 | 4.37E-15 | 4.21E-15 |
| 110 | 1.56E-13 | 1.78E-07 | >1 | 3.06E-14 | 3.01E-14 |
| 210 | 6.21E-13 | >1 | | 5.68E-14 | >1 |
| 310 | 3.28E-13 | | | 8.49E-14 | |
| 410 | 9.04E-13 | | | 2.93E-13 | |
| 510 | 9.51E-13 | | | >1 | |
| 610 | 3.45E-12 | | | | |
| 710 | 5.54E-12 | | | | |
| 810 | 1.56E-11 | | | | |
| 910 | 1.21E-11 | | | | |
| 1010 | 6.71E-11 | | | | |
| 1110 | 1.84E-11 | | | | |
| | … | | | | |
| 1910 | 1.58E-09 | | | | |
| 2010 | 2.90E-09 | | | | |

For $\rho = 1$ both algorithms work equally excellently, for other radii still well.

### 6.2 Disk

Roots and arguments are sampled from the random variables $d_k = \delta\rho e^{i(k+\delta)/n}$ in the complex disk with radius $\rho$.

| | P: Disk | | | | |
|---|---|---|---|---|---|
| n | 1 | 1.2 | 1.4 | 0.5 | 0.1 |
| 10 | 4.15E-13 | 1.30E-13 | 8.38E-13 | 6.75E-13 | 1.49E-14 |
| 20 | 5.03E-11 | 1.43E-12 | 2.13E-12 | 1.48E-12 | 1.95E-10 |
| 30  | 2.60E-08 | 4.87E-10 | 1.17E-11 | 2.32E-08 | 1.27E-10 |
| 40  | 3.31E-06 | 9.10E-07 | 3.08E-10 | 7.69E-11 | 1.07E-08 |
| 50  | 1.65E-09 | 1.53E-07 | 2.54E-06 | 5.53E-08 | 3.92E-07 |
| 60  | 7.34E-08 | 1.95E-06 | 5.35E-06 | 2.98E-07 | 4.45E-07 |
| 70  | 1.02E-07 | 1.20E-06 | 2.66E-04 | 1.43E-08 | 3.71E-07 |
| 80  | 5.91E-06 | 8.57E-04 | 1.19E-05 | 1.56E-04 | 3.66E-05 |
| 90  | >1       | >1       | 2.67E-04 | 2.24E-06 | 1.56E-04 |
| 100 |          |          | 1.08E-02 | >1       | >1       |
| 110 |          |          | >1       |         |         |

| R+:Disk | | | | | |
|---|---|---|---|---|---|
| n | 1 | 1.2 | 1.4 | 0.5 | 0.1 |
| 10  | 4.35E-13 | 1.08E-13 | 8.38E-13 | 6.22E-13 | 1.33E-14 |
| 20  | 5.40E-11 | 1.51E-12 | 2.13E-12 | 1.21E-12 | 1.78E-10 |
| 30  | 2.32E-08 | 5.03E-10 | 1.17E-11 | 2.10E-08 | 1.32E-10 |
| 40  | 3.28E-06 | 1.07E-06 | 3.08E-10 | 9.66E-11 | 1.02E-08 |
| 50  | 1.39E-09 | 1.49E-07 | 2.54E-06 | 5.61E-08 | 3.54E-07 |
| 60  | 7.45E-08 | 2.80E-06 | 5.35E-06 | 2.62E-07 | 4.94E-07 |
| 70  | 1.01E-07 | 1.35E-06 | 2.66E-04 | 1.52E-08 | 3.82E-07 |
| 80  | 5.74E-06 | 9.11E-04 | 1.19E-05 | 1.53E-04 | 3.83E-05 |
| 90  | >1       | >1       | 2.67E-04 | 2.21E-06 | 1.68E-04 |
| 100 |          |          | 1.08E-02 | >1       | >1       |
| 110 |          |          | >1       |         |         |

Both algorithms have almost the same performance.

## 6.3 Line

Roots and arguments are sampled from the random variables $d_k = \delta(-1 + 2(k+\delta)/n)$.

| P::Line | | | | | |
|---|---|---|---|---|---|
| n | 1 | 1.2 | 1.4 | 0.5 | 0.1 |
| 10  | 7.17E-16 | 2.12E-15 | 4.89E-15 | 2.47E-15 | 2.60E-15 |
| 20  | 6.76E-16 | 1.18E-14 | 7.13E-14 | 7.58E-15 | 9.98E-15 |
| 30  | 1.07E-15 | 1.21E-13 | 3.07E-12 | 3.18E-14 | 5.40E-14 |
| 40  | 1.06E-15 | 4.41E-13 | 1.86E-11 | 1.84E-14 | 3.87E-14 |
| 50  | 1.61E-15 | 4.07E-12 | 8.28E-10 | 1.33E-13 | 4.30E-13 |
| 60  | 1.40E-15 | 1.19E-11 | 6.51E-09 | 5.18E-13 | 2.35E-12 |
| 70  | 1.48E-15 | 3.50E-11 | 5.18E-08 | 1.08E-12 | 6.34E-12 |
| 80  | 1.66E-15 | 7.87E-11 | 2.44E-07 | 3.21E-12 | 2.76E-11 |
| 90  | 1.79E-15 | 4.95E-10 | 6.87E-06 | 1.08E-11 | 1.37E-10 |
| 100 | 1.98E-15 | 2.88E-09 | 1.15E-04 | 3.61E-11 | 6.77E-10 |
| 110 | 1.75E-15 | 7.80E-09 | 9.53E-04 | 1.10E-10 | NaN |
| 120 | 1.83E-15 | 5.04E-08 | 1.76E-02 | 2.92E-10 | |
| 130 | 2.14E-15 | 1.38E-07 | 1.05E-01 | 6.83E-10 | |
| 140 | 1.58E-15 | 6.45E-07 | 9.12E-01 | 1.42E-09 | |
| 150 | 1.79E-15 | 5.06E-06 | >1       | 2.60E-09 | |
| 160 | 4.62E-15 | 1.18E-05 |         | 3.78E-09 | |
| 170 | 1.95E-15 | 4.78E-05 |         | 4.00E-09 | |
| 180 | 3.33E-15 | 3.68E-04 |         | 1.12E-08 | |
| 190 | 3.88E-15 | 8.39E-04 |         | 5.56E-08 | |
| 200 | 2.02E-15 | 3.26E-03 |         | 2.10E-07 | |
| 210 | 2.64E-15 | 2.29E-02 |         | 6.63E-07 | |
| 220 | 3.12E-15 | 5.35E-02 |         | NaN      | |
| 230 | 3.09E-15 | 1.71E-01 |         |          | |
| …   | …        | …        |         |          | |
| 360 | 2.31E-15 |          |         |          | |
| 370 | NaN      |          |         |          | |

# FFT-based Computation of Polynomial Coefficients and Related Tasks

| | | | R+:Line | | |
|---|---|---|---|---|---|
| n | 1 | 1.2 | 1.4 | 0.5 | 0.1 |
| 10 | 5.79E-15 | 3.36E-14 | 1.22E-13 | 2.97E-15 | 3.11E-15 |
| 20 | 1.67E-14 | 4.66E-13 | 6.23E-12 | 8.30E-15 | 1.06E-14 |
| 30 | 3.71E-14 | 2.57E-12 | 1.23E-10 | 3.28E-14 | 5.44E-14 |
| 40 | 8.78E-15 | 1.12E-11 | 5.30E-09 | 1.78E-14 | 3.79E-14 |
| 50 | 2.83E-14 | 1.43E-10 | 1.86E-07 | 1.32E-13 | 4.30E-13 |
| 60 | 4.12E-14 | 7.86E-10 | 8.22E-06 | 5.17E-13 | 2.35E-12 |
| 70 | 3.27E-14 | 4.91E-09 | 2.55E-04 | 1.08E-12 | 6.34E-12 |
| 80 | 4.39E-14 | 5.31E-08 | 9.38E-03 | 3.21E-12 | 2.76E-11 |
| 90 | 6.79E-14 | 3.47E-07 | 2.08E-01 | 1.08E-11 | 1.37E-10 |
| 100 | 9.19E-14 | 2.21E-06 | >1 | 3.61E-11 | 6.77E-10 |
| 110 | 1.16E-13 | 9.77E-06 | | 1.10E-10 | 2.93E-09 |
| 120 | 1.30E-13 | 2.84E-05 | | 2.92E-10 | NaN |
| 130 | 1.28E-13 | 8.92E-05 | | 6.83E-10 | |
| 140 | 9.78E-14 | 2.11E-04 | | 1.42E-09 | |
| 150 | 6.92E-14 | 5.88E-03 | | 2.60E-09 | |
| 160 | 4.47E-14 | 4.08E-02 | | 3.78E-09 | |
| 170 | 1.40E-14 | 3.81E-01 | | 4.00E-09 | |
| 180 | 1.95E-14 | >1 | | 1.12E-08 | |
| 190 | 3.93E-14 | | | 5.56E-08 | |
| 200 | 6.93E-14 | | | 2.10E-07 | |
| 210 | 8.49E-14 | | | 6.63E-07 | |
| 220 | 8.98E-14 | | | 4.04E-06 | |
| 230 | 9.37E-14 | | | NaN | |
| … | … | | | | |
| 360 | 6.63E-07 | | | | |
| 370 | 4.04E-06 | | | | |
| 380 | NaN | | | | |

# 7 Conclusion

Algorithms P and R+ have been numerically evaluated on problems A with unit roots as well as on problems F and H using roots and data on circles in the complex plane, inside disks in the complex plane and on the real line. The diameter of the point set has been varied between $0.1$ and $1.4$.

Algorithm P is clearly superior in the case of roots and data on circles, yielding precise results for problem size up to and beyond $n = 2000$. This finding includes problem H. Therefore algorithm PP is excellently suited to invert large Vandermondians and to find the coefficients of polynomials with roots on the unit circle, or circles with radii in the above range.

For roots and data located inside a disk algorithm A retains its performance and its superiority over algorithm R+ on problem F. On problem H both algorithms have similar but a much poorer performance, R+'s slightly better than P's, and become unstable at problem sizes below $n = 100$

With roots and data on the real line the performance of both algorithms decreases further, unstability sets in at $n \leq 40$. Still R+ performs slightly better than P.